\documentclass[12pt,a4paper]{article}
\usepackage{amsmath,amsxtra,amsfonts,amssymb,latexsym,amsthm,amscd,verbatim,amstext}
\usepackage[mathscr]{eucal}

\newcommand{\dif}{\mathrm{d}}

\input cyracc.def

\newfont{\tricyr}{wncyr10 at 12pt}
\newfont{\tricyi}{wncyi10 at 12pt}
\newfont{\tricyb}{wncyb10 at 12pt}
\newfont{\Tricyr}{wncyr10 at 13.6pt}
\newfont{\Tricyi}{wncyi10 at 13.6pt}
\newfont{\Tricyb}{wncyb10 at 13.6pt}
\newfont{\tricmr}{cmr10 at 13.6pt}
\newfont{\tricmi}{cmti10 at 13.6pt}
\newfont{\tricmb}{cmb10 at 13.6pt}

\theoremstyle{plain}
\newtheorem{Th}{Theorem}

\newtheorem{Lem}{Lemma}

\theoremstyle{definition}

\newtheorem{Not}{Remark}

\begin{document}
\centerline {{\bf On the initial value problem for the Navier-Stokes equations with}} 
\centerline {{\bf the initial datum in critical Sobolev and Besov spaces}}

\vskip 0.7cm
\begin{center}
D. Q. Khai, N. M. Tri
\end{center}
\begin{center}
Institute of Mathematics, Vietnam Academy of Science and Technology\\
18 Hoang Quoc Viet, 10307 Cau Giay, Hanoi, Vietnam
\end{center}
\vskip 0.3cm 
{\bf Abstract}: The existence of local unique mild 
solutions to the Navier-Stokes
equations in the whole space with an initial tempered distribution 
datum in critical homogeneous or inhomogeneous Sobolev spaces is shown. 
Especially, the case when the integral-exponent is less than 2 is investigated. 
The global existence is also obtained for the initial datum in critical homogeneous  
Sobolev spaces with a norm small enough in  suitable critical  Besov spaces. 
The key lemma is to establish the bilinear estimates in these spaces, 
due to the point-wise decay of the kernel of the heat semigroup.



\footnotetext[3]{{\it e-mail address}: khaitoantin@gmail.com triminh@math.ac.vn}
\vskip 1.0cm 
\centerline{\S1. Introduction} 
\vskip 0.5cm 
We consider the Navier-Stokes equations (NSE) in $d$ 
dimensions in special setting of a viscous, homogeneous, 
incompressible fluid which fills the entire space and is 
not submitted to external forces. Thus, the equations 
we consider are the system: 
\begin{align} 
\left\{\begin{array}{ll} \partial _tu  = \Delta u - 
\nabla .(u \otimes u) - \nabla p , & \\ 
{\rm div}(u) = 0, & \\
u(0, x) = u_0, 
\end{array}\right . \notag
\end{align}
which is a condensed writing for
\begin{align} 
\left\{\begin{array}{ll} 1 \leq k \leq d, \ \  \partial _tu_k  
= \Delta u_k - \sum_{l =1}^{d}\partial_l(u_lu_k) - \partial_kp , & \\ 
\sum_{l =1}^{d}\partial_lu_l = 0, & \\
1 \leq k \leq d, \ \ u_k(0, x) = u_{0k} .
\end{array}\right . \notag
\end{align}
The unknown quantities are the velocity $u(t, x)=(u_1(t, x),\dots,u_d(t, x))$ 
of the fluid element at time $t$ and position $x$ and the pressure $p(t, x)$.\\
A translation invariant Banach space of tempered distributions $\mathcal{E}$ 
is called a critical space for NSE if its norm is invariant under the action 
of the scaling $f(.) \longrightarrow \lambda f(\lambda .)$. One can take, 
for example, $\mathcal{E} = L^d(\mathbb{R}^d)$ 
or the smaller space $\mathcal{E} = \dot{H}^{\frac{d}{2} - 1}(\mathbb{R}^d)$. 
In fact, one has the "chain of critical spaces" given by the continuous embeddings
\begin{equation}\label{Critical spaces}
\dot{H}^{\frac{d}{2} - 1}(\mathbb{R}^d) \hookrightarrow L^d(\mathbb{R}^d) \hookrightarrow \dot{B}^{\frac{d}{q}-1, \infty}_{q}(\mathbb{R}^d)_{(d\leq q < \infty)} \hookrightarrow BMO^{-1}(\mathbb{R}^d) \hookrightarrow \dot{B}^{- 1, \infty}_{\infty}(\mathbb{R}^d).
\end{equation}
It is remarkable feature that NSE are well-posed in the sense of Hadarmard 
(existence, uniqueness and continuous dependence on data) when the initial datum 
is divergence-free and belong to the critical function spaces 
(except $\dot{B}^{- 1, \infty}_{\infty}$) listed in $\eqref{Critical spaces}$ 
(see \cite{M. Cannone 1995} for $\dot{H}^{\frac{d}{2} - 1}(\mathbb{R}^d)$, 
$L^d(\mathbb{R}^d)$, and $\dot{B}^{\frac{d}{q}-1, \infty}_{q}(\mathbb{R}^d)$,  
see \cite{H. Koch 2001} for $BMO^{-1}(\mathbb{R}^d)$. The recent 
ill-posedness result  for 
$\dot{B}^{- 1, \infty}_{\infty}(\mathbb{R}^d)$) with $d \geq 3$ was established 
in \cite{J. Bourgain 2008}. However, the ill-posedness in $\dot{B}^{- 1, \infty}_{\infty}(\mathbb{R}^d)$  is still open when $d = 2$.\\
In the 1960s, mild solutions were first  constructed by Kato and 
Fujita (\cite{T. Kato 1962}, \cite{H. Fujita 1964}) that are continuous 
in time and take values in the Sobolev space $H^s(\mathbb{R}^d),\linebreak
(s \geq \frac{d}{2} - 1)$, say $u \in C([0, T]; H^s(\mathbb{R}^d))$. 
In 1992, a modern treatment for mild solutions in 
$H^s(\mathbb{R}^d), (s \geq \frac{d}{2} - 1)$ 
was given by Chemin \cite{J. M. Chemin 1992}. In 1995, using the 
simplified version of the bilinear operator, Cannone proved 
the existence of mild solutions in $\dot{H}^s(\mathbb{R}^d), 
(s \geq \frac{d}{2} - 1)$, see \cite{M. Cannone 1995}. 
Results on the existence of mild solutions with value in 
$L^q(\mathbb{R}^d), (q > d)$ were established in  the papers 
of  Fabes, Jones and Rivi\`{e}re \cite{E. Fabes 1972a} 
and of Giga \cite{Y. Giga 1986a}. Concerning the initial 
datum in the space $L^\infty(\mathbb R^d)$, the existence of a mild solution was 
obtained by Cannone and Meyer in (\cite{M. Cannone 1995}, 
\cite{M. Cannone. Y. Meyer 1995}). Moreover, in (\cite{M. Cannone 1995}, 
\cite{M. Cannone. Y. Meyer 1995}), they also obtained theorems 
on the existence of mild solutions with value in Morrey-Campanato space 
$M^q_2(\mathbb{R}^d), (q > d)$ and Sobolev space 
$H^s_q(\mathbb{R}^d), (q < d, \frac{1}{q} 
- \frac{s}{d} < \frac{1}{d})$, and in general in the case of a so-called 
well-suited space $\mathcal{W}$ for NSE. NSE in the  
Morrey-Campanato spaces were also treated by Kato  \cite{T. Kato 1992},  
Taylor \cite{M.E. Taylor 1992}, Kozono and Yamazaki 
\cite{H. Kozono: 1994a}.\\
In 1981, Weissler \cite{F.B. Weissler:1981} gave the first existence 
result of mild solutions in the half space $L^3(\mathbb{R}^3_+)$. 
Then Giga and Miyakawa \cite{Y. Giga 1985} generalized the
result to $L^3(\Omega)$, where $\Omega$ is an open bounded domain 
in $\mathbb{R}^3$. Finally, in 1984, Kato \cite{T. Kato 1984} obtained, 
by means of a purely analytical tool (involving only the H$\ddot{\text{o}}$lder and 
Young inequalities and without using any estimate of fractional powers 
of the Stokes operator), an existence theorem in the whole space 
$L^3(\mathbb{R}^3)$. In (\cite{M. Cannone 1995}, 
\cite{M. Cannone 1997}, \cite{M. Cannone 1999}), 
Cannone showed how to simplify Kato's proof. The idea is 
to take advantage of the structure of the bilinear operator in its scalar form. 
In particular, the divergence $\nabla$ and heat $e^{t\Delta}$ operators can 
be treated as a single convolution operator. In 1994, Kato and Ponce \cite{T. Kato 1994} 
showed that NSE are well-posed when the initial datum belongs to the 
homogeneous Sobolev spaces $\dot{H}^{\frac{d}{q} - 1}_q(\mathbb{R}^d), 
(d \leq q < \infty)$. 
Recently, the authors of this article have 
considered NSE in mixed-norm Sobolev-Lorentz spaces, 
see \cite{N. M. Tri: Tri2014a}. In  \cite{N. M. Tri: Tri2014???}, 
we showed that NSE are well-posed when  the initial datum belongs 
to the Sobolev spaces $\dot{H}^s_q(\mathbb{R}^d)$ with non-positive-regular 
indexes $ (q \geq d, \frac{d}{q}-1 \leq s \leq 0)$. In \cite{N. M. Tri: Tri2014??}, 
we showed that the bilinear operator
\begin{equation}\label{B}
B(u, v)(t) = \int_{0}^{t} e^{(t-\tau ) \Delta} \mathbb{P}
 \nabla.\big(u(\tau,.)\otimes v(\tau,.)\big) \dif\tau
\end{equation}
is bicontinuous in $L^\infty([0, T]; \dot{H}^s_q(\mathbb{R}^d))$ 
with super-critical, non-negative-regular indexes $(0 \leq s \leq d - 1, q > 1,\ {\rm and}\ \frac{s}{d}<\frac{1}{q}<{\rm min}\Big\{\frac{s+1}{d},\frac{s+d}{2d}\Big\})$, 
and we obtain the inequality
\begin{gather*}
\big\|B(u, v)\big\|_{L^\infty([0, T]; \dot{H}^s_q)} \leq
 C_{s,q,d}T^{\frac{1}{2}(1 + s - \frac{d}{q})}\big\|u\big\|_{L^\infty([0, T]; \dot{H}^s_q)}
\big\|v\big\|_{L^\infty([0, T]; \dot{H}^s_q)}. 
\end{gather*}
In this case existence and uniqueness theorems of local mild 
solutions can therefore be easily deduced.\\
In this paper, first, for $d \geq 3, s \geq 0,  p > 1, {\rm and}\ r > 2$ 
be such that $\frac{s}{d} < \frac{1}{p} < \frac{1}{2} +\frac{s}{2d}
\ {\rm and}\  \frac{2}{r} + \frac{d}{p} - s  \leq 1$, we investigate 
mild solutions to NSE in the spaces $L^r\big([0, T]; 
\dot{H}^s_p(\mathbb{R}^d)\big)$. 
We obtain the existence of local mild solutions with arbitrary 
initial tempered 
distribution datum in the Besov spaces $B_p^{s - \frac{2}{r},r}$. 
In the case of critical indexes $\frac{2}{r} - s + \frac{d}{p} = 1$, 
we obtain the existence of global mild solutions when the norm 
of the initial tempered 
distribution datum in the Besov space $\dot{B}_p^{s -\frac{2}{r},r}$ 
is small enough. The particular case of the above result, when 
$s=0$, was presented in the book by Lemarie-Rieusset 
\cite{P. G. Lemarie-Rieusset 2002}. We also note that the 
Cauchy problem for an incompressible 
magneto-hydrodynamics system with positive viscosity and magnetic resistivity, 
in the framework of the Besov spaces was considered in \cite{C. Miao 2002}. \\
 Next, we present two different algorithms 
for constructing mild solutions in $C([0, T];\dot{H}^{\frac{d}{q}-1}_q(\mathbb{R}^d))$ 
or $C([0, T];H^{\frac{d}{q}-1}_q(\mathbb{R}^d))$  
to the Cauchy problem for the Navier-Stokes equations when  
the initial datum belongs to the Sobolev spaces 
$\dot{H}^{\frac{d}{q}-1}_q(\mathbb{R}^d)\  
\big({\rm or}\ H^{\frac{d}{q}-1}_q(\mathbb{R}^d)\big)$. 
We use the first algorithm to consider the case when 
the initial datum belongs to $\dot{H}^{\frac{d}{q} - 1}_q(\mathbb{R}^d)$ or
$H^{\frac{d}{q} - 1}_q(\mathbb{R}^d)$ 
with $3 \leq d \leq 4$ and $2 \leq q \leq d$. Our results, when  $q = d$, 
are a generalization the ones  obtained in \cite{P. G. Lemarie-Rieusset 2002}.  
With the second algorithm, we can treat the case when 
the initial datum belongs to the critical spaces 
$\dot{H}^{\frac{d}{q} - 1}_q(\mathbb{R}^d)$ with $d \geq3$ and $1 < q \leq d$.
The cases $q = 2$ and $q = d$ were considered by many authors, 
see (\cite{M. Cannone 1995}, \cite{M. Cannone 1999},  \cite{J. M. Chemin 1992},  \cite{Hongjie Dong 2007}, \cite{H. Fujita 1964}, \cite{T. Kato 1962}, \cite{T. Kato 1984}, \cite{P. G. Lemarie-Rieusset 2002}, \cite{F. Planchon 1996}). 
A part of our results  in the case when 
$2 < q < d$ can also be obtained by using the interpolation method  
of the results between the spaces $\dot H^{\frac{d}{2}}$ and $L^d$. \\
So we will concentrate our efforts on the case $1 < q < 2$. 
To obtain the existence theorem in $C([0, T];\dot{H}^{\frac{d}{q}-1}_q(\mathbb{R}^d))$, 
we need to establish the continuity of the bilinear operator $B$ from 
$$
L^{2q}\Big([0, T]; \dot{H}^{\frac{d + 2 - 2q}{q} }_{\frac{dq}{d + 1 - q}}\Big) 
\times L^{2q}\Big([0, T]; \dot{H}^{\frac{d + 2 
- 2q}{q} }_{\frac{dq}{d + 1 - q}}\Big) \ {\rm to }\ \ C([0, T];
\dot{H}^{\frac{d}{q}-1}_q(\mathbb{R}^d)),
$$
and establishes the continuity of the bilinear operator $B$  from
$
L^r([0, T]; H^s_p) \times L^r([0, T]; H^s_p)
$
into 
$
L^r([0, T]; H^s_p).
$  
In order to evaluate the norm of the bilinear operator $B$ in these spaces 
we use Lemma \ref{Them} which estimates the point-wise 
product of two functions in  $\dot{H}^s_q(\mathbb{R}^d)$.

The paper is organized as follows. In Section 2 we recall some 
embedding theorems in the Triebel and Besov spaces and auxiliary lemmas. 
In Section 3 we present the main results of the paper. \\
In the sequence, for a space of functions defined on $\mathbb R^d$, 
say $E(\mathbb R^d)$, we will abbreviate it as $E$. 
\vskip 1cm
\centerline{\S2. Some imbedding theorems}
\vskip 0.7cm
In this paper we use the definition of the Besov space  $B^{s,p}_q$, 
the Triebel space $F^{s,p}_q$, and their homogeneous space 
$\dot B^{s,p}_q$ and $\dot F^{s,p}_q$  in \cite{G. Bourdaud 1988, 
G. Bourdaud 1993, M. Frazier 1991, J. Peetre 1976}. A known property 
of these spaces is the Riesz potential $\dot \Lambda^s = (-\Delta)^{s/2}$ 
which is an isomorphism from $\dot B^{s_0,p}_q$ onto $\dot B^{s_0-s,p}_q$ 
and from $\dot F^{s_0,p}_q$ to $\dot F^{s_0-s,p}_q$, see \cite{Bjorn Jawerth: 1977}.\\
Let $1 < q < \infty$ and $s < d/q$, we define the homogeneous 
Sobolev space $\dot{H}^s_q$ as the closure of the space
$S_0 = \big\{f \in \mathcal{S}: \ 0 \notin {\rm Supp }\hat f\big \}$ 
in the norm $\|f\|_{\dot H^s_q} = \|\dot \Lambda^s f\|_q$. 
Let us recall the following lemmas.
\begin{Lem}\label{lem.2.1.4} 
Let $1 \leq p, q \leq \infty$ and $s \in \mathbb{R}$.\\
{\rm (a)} If $ s < 1$ then the two quantities
\begin{gather}
\ \Big(\int_0^\infty\big(t^{-\frac{s}{2}}\big\|e^{t\Delta}
t^{\frac{1}{2}}\dot{\Lambda}f\big\|_q\big)^p\frac{\dif t}{t} \Big)^{1/p} 
\ and \ \big\|f\big\|_{\dot{B}_{q}^{s, p}}
 \ are \ equivalent. \notag
\end{gather}
{\rm (b)} If $ s < 0$ then the two quantities 
\begin{gather}
\ \Big(\int_0^\infty\big(t^{-\frac{s}{2}}\big\|e^{t\Delta}f\big\|_q\big)^p
\frac{\dif t}{t} \Big)^{1/p} \ and \ \big\|f\big\|_{\dot{B}_{q}^{s, p}} 
\ are \ equivalent. \notag
\end{gather}
\end{Lem}
\textbf{Proof}: See (\cite{S. Friedlander 2004}, Proposition 1, p. 181
 and Proposition 3, p. 182), or see (\cite{P. G. Lemarie-Rieusset 2002}, 
Theorem 5.4, p. 45).\qed \\
The following lemma is a generalization of the above lemma.
\begin{Lem}\label{lem.2.1.4'} 
Let $1 \leq p, q \leq \infty,\ \alpha \geq 0$, and $ s < \alpha$. Then the two quantities
\begin{gather*}
\Big(\int_0^\infty(t^{-\frac{s}{2}}\big\|e^{t\Delta}t^{\frac{\alpha}{2}}
\dot \Lambda^\alpha f\big\|_{L^q})^p\frac{{\rm d}t}{t}\Big)^{\frac{1}{p}}
\ and \ \big\|f\big\|_{\dot{B}_{q}^{s, p}} \ are \ equivalent. 
\end{gather*}
\end{Lem}
\textbf{Proof}: Note that $\dot{\Lambda}^{s_0}$ is an isomorphism 
from $\dot{B}^{s,p}_q$ to $\dot{B}^{s-s_0,p}_q$, then we can easily 
prove the lemma.   \qed
\begin{Lem}\label{lem.2.2.3}
For $1 \leq p, q,r \leq \infty$  and $s \in \mathbb{R}$, 
we have the following embedding mappings.\\
{\rm (a)} If $1 < q \leq 2$ then 
$$
\dot{B}_{q}^{s, q} \hookrightarrow \dot{H}_{q}^s 
\hookrightarrow \dot{B}_{q}^{s, 2}, \ B_{q}^{s, q} 
\hookrightarrow H_{q}^s \hookrightarrow B_{q}^{s, 2}.
$$
{\rm (b)} If $2 \leq  q < \infty$ then 
$$
\dot{B}_{q}^{s, 2} \hookrightarrow \dot{H}_{q}^s 
\hookrightarrow \dot{B}_{q}^{s, q}, \ B_{q}^{s, 2} 
\hookrightarrow H_{q}^s \hookrightarrow B_{q}^{s, q}.
$$
{\rm (c)} If $1 \leq p_1 < p_2 \leq \infty$  then
$$
\dot{B}_{q}^{s, p_1} \hookrightarrow  \dot{B}_{q}^{s, p_2}, B_{q}^{s, p_1} 
\hookrightarrow  B_{q}^{s, p_2}, \ \dot{F}_{q}^{s, p_1} 
\hookrightarrow  \dot{F}_{q}^{s, p_2}, F_{q}^{s, p_1} 
\hookrightarrow  F_{q}^{s, p_2}.
$$
{\rm (d)} If $s_1 > s_2, \ 1 \leq q_1, \ q_2\leq \infty$,  
and $s_1 - \frac{d}{q_1} = s_2 - \frac{d}{q_2}$
then
$$
\dot{B}_{q_1}^{s_1, p} \hookrightarrow \dot{B}_{q_2}^{s_2, p}, B_{q_1}^{s_1, p} 
\hookrightarrow B_{q_2}^{s_2, p}, \dot{F}_{q_1}^{s_1, p} 
\hookrightarrow \dot{F}_{q_2}^{s_2, r}, \ F_{q_1}^{s_1, p} 
\hookrightarrow F_{q_2}^{s_2, r}.
$$
{\rm (e)} If $p \leq q$  then
$$
B_{q}^{s, p} \hookrightarrow  F_{q}^{s, p}, 
\ \dot{B}_{q}^{s, p} \hookrightarrow  \dot{F}_{q}^{s, p}. 
$$
{\rm (f)} If $q \leq p$  then
$$
F_{q}^{s, p} \hookrightarrow  B_{q}^{s, p}, \ \dot{F}_{q}^{s, p} 
\hookrightarrow  \dot{B}_{q}^{s, p}. 
$$
{\rm (g)}
$$
F_{q}^{s, q} =  B_{q}^{s, q}, \ \dot{F}_{q}^{s, q} =  \dot{B}_{q}^{s, q}. 
$$
{\rm (h)} If $1 < q < \infty$
$$
H_q^s = F_q^{s, 2} , \ \dot{H}_q^s =  \dot{F}_q^{s, 2}. 
$$
\end{Lem}
\textbf{Proof}: For the proof of (a) and (b) see Theorem 6.4.4 
(\cite{J. Bergh: 1976}, p. 152). For the proof of (c) see \cite{A. R. Adams:1975a} 
and \cite{J. Bergh: 1976}. For the proof of (d) see Theorem 6.5.1 
(\cite{J. Bergh: 1976}, p. 153) and \cite{Bjorn Jawerth: 1977}.  
For the proof of (e), (f), (g), and (h) see \cite{A. R. Adams:1975a} 
and \cite{Bjorn Jawerth: 1977}. \qed
\begin{Lem}\label{lem.2.2.6}
Let $p \geq 1\ and\ s \in \mathbb{R}$. Then the following statements hold\\
{\rm (1)} Assume that $u_0 \in H^s_p$. Then\\
$e^{t\Delta}u_0 \in L^\infty([0, \infty); H^s_p) 
\ \text{and}\ \big\|e^{t\Delta}u_0\big\|_{L^\infty([0, \infty); H^s_p)} 
\leq \big\|u_0\big\|_{H^s_p}.$\\
{\rm (2)} Assume that $u_0 \in \dot{H}^s_p$. Then \\
$e^{t\Delta}u_0 \in L^\infty([0, \infty); \dot{H}^s_p) 
\ \text{and}\ \big\|e^{t\Delta}u_0\big\|_{L^\infty([0, \infty); \dot{H}^s_p)} 
\leq\big\|u_0\big\|_{\dot{H}^s_p}.$
\end{Lem}
\textbf{Proof}: (1) We have 
\begin{gather*}
\big\|e^{t\Delta}u_0\big\|_{H^s_p} = \big\|e^{t\Delta}(Id-\Delta)^{s/2}
u_0\big\|_{L^p} = \notag \\ 
\frac{1}{(4\pi t)^{d/2}} \Big\|\int_{\mathbb R^d}
e^{\frac{-|\xi|^2}{4t}}\big((Id-\Delta)^{s/2}u_0\big)(\ . -\xi)\dif\xi\Big\|_{L^p} \\ \leq
 \frac{1}{(4\pi t)^{d/2}}\int_{\mathbb R^d}
e^{\frac{-|\xi|^2}{4t}}\big\|\big((Id-\Delta)^{s/2}u_0\big)(\ . -\xi)\big\|_{L^p}\dif\xi \\
  = \frac{1}{(4\pi t)^{d/2}}\int_{\mathbb R^d}
e^{\frac{-|\xi|^2}{4t}}\big\|u_0\big\|_{H^s_p}\dif\xi 
= \big\|u_0\big\|_{H^s_p}, \ t \geq 0 . 
\end{gather*}
(2) The proof of (2) is similar to the proof of (1). \qed 
\begin{Th}\label{Th.2.1}
Let $E$ be an Banach space, and let $B: E \times E \rightarrow  E$ 
be a continuous bilinear form such that there exists $\eta > 0$ so that
$$
\|B(x, y)\| \leq \eta \|x\| \|y\|,
$$
for all x and y in $E$. Then for any fixed $y \in E$ 
such that $\|y\| \leq \frac{1}{4\eta}$, the equation $x = y - B(x,x)$ 
has a unique solution  $\overline{x} \in E$ satisfying $ \|\overline{x}\| \leq \frac{1}{2\eta}$.
\end{Th}
\textbf{Proof}: See Theorem 22.4 (\cite{P. G. Lemarie-Rieusset 2002}, p. 227).\qed\\

The following lemmas, in which we estimate the point-wise product of two 
functions in  $\dot{H}^s_p(\mathbb{R}^d)$ is more general than the H\"{o}lder inequality. 
In the case when $s = 0, p \geq 2,$ we get back the usual  H\"{o}lder inequality. 

\begin{Lem}\label{lem2.1}
Assume that 
$$
1 < p, q < d\  and\  \frac{1}{p} + \frac{1}{q} < 1 +  \frac{1}{d}.
$$
Then there exists a constant $C$ independent of $u,v$ such that 
the following inequality holds
$$
\big\|uv\big\|_{\dot{H}^1_r}  \leq C\big\|u\big\|_{\dot{H}^1_p}\big\|v\big\|_{\dot{H}^1_q}, 
\ \forall  u \in \dot{H}^1_p, v \in \dot{H}^1_q,
$$
where $\frac{1}{r} = \frac{1}{p} + \frac{1}{q} - \frac{1}{d}$. In the subsequence 
the above kinds of conclusions will be shorten as 
$$
\big\|uv\big\|_{\dot{H}^1_r}  \lesssim \big\|u\big\|_{\dot{H}^1_p}
\big\|v\big\|_{\dot{H}^1_q}.
$$
\end{Lem}
\textbf{Proof}: 
By applying the Leibniz formula for the derivatives of a product 
of two functions, we have 
\begin{gather*}
 \big\|uv\big\|_{\dot{H}^1_r} \simeq \sum_{|\alpha| = 1}
\big\|\partial^{\alpha}(uv)\big\|_{L^r} \leq \sum_{|\alpha| = 1}
\big\|(\partial^{\alpha}u)v\big\|_{L^r} + \sum_{|\alpha| = 1}
\big\|u(\partial^{\alpha}v)\big\|_{L^r}.
\end{gather*}
From the H\"{o}lder and Sobolev inequalities it follows that  
\begin{gather*}
\sum_{|\alpha| = 1}\big\|(\partial^{\alpha}u)v\big\|_{L^r} 
\leq \sum_{|\alpha| = 1}\big\|\partial^{\alpha}u\big\|_{L^p}
\big\|v\big\|_{L^{q_1}} \lesssim \big\|u\big\|_{\dot{H}^1_p}
\big\|v\big\|_{\dot{H}^1_q},
\end{gather*}
where
$$
\frac{1}{q_1} = \frac{1}{q} - \frac{1}{d}.
$$
Similar to the above proof, we have
\begin{gather*}
\sum_{|\alpha| = 1}\big\|u(\partial^{\alpha}v)\big\|_{L^r}  
\lesssim \big\|u\big\|_{\dot{H}^1_p}\big\|v\big\|_{\dot{H}^1_q}.
\end{gather*}
This gives the desired result
$$
\big\|uv\big\|_{\dot{H}^1_r} \lesssim \big\|u\big\|_{\dot{H}^1_p}
\big\|v\big\|_{\dot{H}^1_q}.\qed
$$
\begin{Lem}\label{lem2.1'}
Assume that 
\begin{equation}\label{them1}
0 \leq s \leq 1 , \frac{1}{p} > \frac{s}{d}, \frac{1}{q} > \frac{s}{d},
\  and\  \frac{1}{p} + \frac{1}{q} < 1 +  \frac{s}{d}.
\end{equation}
Then the following inequality holds
$$
\big\|uv\big\|_{\dot{H}^s_r} \lesssim \big\|u\big\|_{\dot{H}^s_p}
\big\|v\big\|_{\dot{H}^s_q}, 
\ \forall  u \in \dot{H}^s_p, v \in \dot{H}^s_q,
$$
where $\frac{1}{r} = \frac{1}{p} + \frac{1}{q} - \frac{s}{d}$.
\end{Lem}
\textbf{Proof}: It is not difficult to show that if $p, q, {\rm and}\ s$ 
satisfy \eqref{them1} then there exist numbers $p_1,p_2,q_1,q_2\in(1,+\infty)$ 
(may be many of them) such that 
\begin{gather*}
\frac{1}{p} = \frac{1 - s}{p_1} + \frac{s}{p_2},  \frac{1}{q} 
= \frac{1 - s}{q_1} + \frac{s}{q_2}, \frac{1}{p_1} + \frac{1}{q_1} < 1, \\ 
 p_2 < d, q_2 < d,\ {\rm and}\ \frac{1}{p_2} + \frac{1}{q_2} < 1 +  \frac{1}{d}.
\end{gather*}
Setting
$$
\frac{1}{r_1} = \frac{1}{p_1} + \frac{1}{q_1}, \frac{1}{r_2} 
= \frac{1}{p_2} + \frac{1}{q_2} - \frac{1}{d}, 
$$
we have
$$
\frac{1}{r} = \frac{1 - s}{r_1} + \frac{s}{r_2}.
$$
Therefore, applying Theorem 6.4.5 (p. 152) of \cite{J. Bergh: 1976} 
(see also \cite{N. Kalton: 2007a} for $\dot H_{p}^{s}$), we get
\begin{gather*}
\dot{H}^s_p = [L^{p_1}, \dot{H}^1_{p_2}]_{s},  \dot{H}^s_q 
= [L^{q_1}, \dot{H}^1_{q_2}]_{s}, \dot{H}^s_r 
= [L^{r_1}, \dot{H}^1_{r_2}]_{s}.  
\end{gather*}
Applying the H\"{o}lder inequality and Lemma \ref{lem2.1} 
in order to obtain 
\begin{gather*}
\big\|uv\big\|_{L^{r_1}} \lesssim \big\|u\big\|_{L^{p_1}}\big\|v\big\|_{L^{q_1}},
\ \forall  u \in L^{p_1}, v \in L^{q_1},\\
\big\|uv\big\|_{\dot{H}^1_{r_2}} \lesssim \big\|u\big\|_{\dot{H}^1_{p_2}}
\big\|v\big\|_{\dot{H}^1_{q_2}}, 
\ \forall  u \in \dot{H}^1_{p_2}, v \in \dot{H}^1_{q_2}.
\end{gather*}
From  Theorem 4.4.1 (p. 96) of \cite{J. Bergh: 1976} we get
$$
\big\|uv\big\|_{{\dot{H}^s_r}} \lesssim \big\|u\big\|_{\dot{H}^s_p}
\big\|v\big\|_{\dot{H}^s_q}.\qed
$$
 \begin{Lem}\label{Them}
Assume that  
\begin{equation}\label{Them1}
0 \leq s < d, \frac{s}{d} < \frac{1}{p}, \frac{s}{d} <  
\frac{1}{q},\ and\ \frac{1}{p} + \frac{1}{q}  <  1 + \frac{s}{d}. 
\end{equation} 
Then we have the inequality 
\begin{gather*}
\big\|uv\big\|_{\dot{H}^s_r} \lesssim \big\|u\big\|_{\dot{H}^s_p}
\big\|v\big\|_{\dot{H}^s_q},
\ \forall  u \in \dot{H}^s_p,  v \in \dot{H}^s_q,
\end{gather*}
where $\frac{1}{r} = \frac{1}{p} + \frac{1}{q} - \frac{s}{d}$.
\end{Lem}
\textbf{Proof}: Denote by $[s]$ the integer part of $s$ and by $\{s\}$ 
the fraction part of $s$. Using formula for the 
derivatives of a product of two functions, we have
\begin{gather*}
\big\|uv\big\|_{\dot{H}^s_r} = \big\|\dot{\Lambda}^s(uv)\big\|_{L^r} 
= \big\|\dot{\Lambda}^{\{s\}}(uv)\big\|_{\dot{H}^{[s]}_r} \simeq \\
\sum_{|\alpha| = [s] }\big\|\partial^{\alpha}\dot{\Lambda}^{\{s\}}(uv)\big\|_{L^r}  = 
\sum_{|\alpha| = [s] }\big\|\dot{\Lambda}^{\{s\}}\partial^{\alpha}(uv)\big\|_{L^r}\\
= \sum_{|\alpha| = [s] }\big\|\partial^{\alpha}(uv)\big\|_{\dot{H}^{\{s\}}_r}
 \lesssim  \sum_{|\gamma| + |\beta| = [s]}
\big\|\partial^{\gamma}u\partial^{\beta}v\big\|_{\dot{H}^{\{s\}}_r}.
\end{gather*}
Set
$$
\frac{1}{\tilde p} = \frac{1}{p} - \frac{s -|\gamma| - \{s\}}{d}, 
\frac{1}{\tilde q} = \frac{1}{q} - \frac{s -|\beta| - \{s\}}{d}.   
$$
Applying Lemma \ref{lem2.1'} and the Sobolev inequality in order to obtain 
\begin{gather*}
\big\|\partial^{\gamma}u\partial^{\beta}v\big\|_{\dot{H}^{\{s\}}_r} 
\lesssim \big\|\partial^{\gamma}u\big\|_{\dot{H}^{\{s\}}_{\tilde p}} 
\big\|\partial^{\beta}v\big\|_{\dot{H}^{\{s\}}_{\tilde q}} \\
\lesssim \big\|u\big\|_{\dot{H}^{|\gamma| + \{s\}}_{\tilde p}} 
\big\|v\big\|_{\dot{H}^{|\beta| + \{s\}}_{\tilde q}} 
\lesssim \big\|u\big\|_{\dot{H}^s_p}\big\|v\big\|_{\dot{H}^s_q}.
\end{gather*}
This gives the desired result
$$
\big\|uv\big\|_{\dot{H}^s_r} \lesssim \big\|u\big\|_{\dot{H}^s_p}
\big\|v\big\|_{\dot{H}^s_q}.\qed
$$
\begin{Not}
Lemmas \ref{lem2.1}, \ref{lem2.1'}, and \ref{Them} 
are still valid when the homogeneous space $\dot{H}^s_p$ 
is replaced by the inhomogeneous space $H^s_p$.  
\end{Not}
\vskip 0.5cm 
\centerline{\S3. The main results}
\vskip 0.5cm 
For $T > 0$, we say that $u$ is a mild solution of NSE on $[0, T]$ 
corresponding to a divergence-free initial data $u_0$ when $u$ 
satisfies the integral equation
$$
u = e^{t\Delta}u_0 - \int_{0}^{t} e^{(t-\tau) \Delta} \mathbb{P} 
\nabla  .\big(u(\tau,.)\otimes u(\tau,.)\big) \dif\tau.
$$
Above we have used the following notation: For a tensor $F = (F_{ij})$ 
we define the vector $\nabla.F$ by $(\nabla.F)_i = \sum_{i = 1}^d\partial_jF_{ij}$ 
and for vectors $u$ and $v$, we define their tensor product $(u \otimes v)_{ij} = u_iv_j$. 
The operator $\mathbb{P}$ is the Helmholtz-Leray projection 
onto the divergence-free fields 
\begin{equation}\label{3.1} 
(\mathbb{P}f)_j =  f_j + \sum_{1 \leq k \leq d} R_jR_kf_k, 
\end{equation} 
where $R_j$ is the Riesz transforms defined on a scalar function $g$ as 
$$
\widehat{R_jg}(\xi) = \frac{i\xi_j}{|\xi|}\hat{g}(\xi).
$$
The heat kernel $e^{t\Delta}$ is defined as 
$$
e^{t\Delta}u(x) = ((4\pi t)^{-d/2}e^{-|.|^2/4t}*u)(x).
$$
If $X$ is a normed space and $u = (u_1, u_2,...,u_d), 
u_i \in X, 1 \leq i \leq d$, then we write 
$$
u \in X, \|u\|_X = \Big(\sum_{i = 1}^d\|u_i\|_X^2\Big)^{1/2}.
$$
\vskip 0.5cm
\textbf{3.1. On the continuity and regularity of the bilinear operator}
\vskip 0.3cm
In this subsection a particular attention will be devoted to the study of 
the bilinear operator $B(u,v)(t)$ defined by \eqref{B}.
\begin{Lem}\label{lem.2.1.2}
Let 
\begin{equation}\label{lem.2.1.2.a}
d \geq 3, \ s \geq 0, \ p > 1, \ r > 2, and\ T>0
\end{equation}
be such that
\begin{equation}\label{lem.2.1.2.b}
\frac{s}{d} < \frac{1}{p} < \frac{1}{2} +\frac{s}{2d}
\ and\  \frac{2}{r} + \frac{d}{p} - s  \leq 1.
\end{equation}
Then the bilinear operator  $B(u, v)(t)$ is continuous from 
$$
L^r([0, T]; H^s_p) \times L^r([0, T]; H^s_p)
$$
into 
$$
L^r([0, T]; H^s_p),
$$
and the following inequality holds
\begin{equation}\label{2.1.1}
\big\|B(u, v)\big\|_{L^r([0, T]; H^s_p)} 
\leq CT^{\frac{1}{2}(1 + s - \frac{2}{r} - \frac{d}{p})}
\big\|u\big\|_{L^r([0, T]; H^s_p)}\big\|v\big\|_{L^r([0, T]; H^s_p)},
\end{equation}
where C is a positive constant independent of T.
\end{Lem}
\textbf{Proof}: We have 
\begin{gather}
\big\|B(u,v)(t)\big\|_{H^s_p} \leq \int_{0}^{t} \Big\|e^{(t- \tau) \Delta}
 \mathbb{P} \nabla .\big(u(\tau,.)\otimes v(\tau,.)\big)\Big\|_{H^s_p} \dif \tau  = \notag \\
 \int_{0}^{t} \Big\|e^{(t- \tau) \Delta} \mathbb{P} \nabla .(Id-\Delta)^{s/2}
\big(u(\tau,.)\otimes v(\tau,.)\big)\Big\|_{L^p}\dif \tau ,\label{2.1.2}
\end{gather}
where the operator $(Id - \Delta)^{\frac{s}{2}}$ is defined via 
the Fourier transform as  
$$
\big((Id - \Delta)^{\frac{s}{2}}g\big)^{\land}(\xi) 
= (1 +|\xi|^2 )^{\frac{s}{2}}\hat{g}(\xi).
$$
We have
\begin{gather}
\Big(e^{(t - \tau)\Delta}\mathbb{P}\nabla .(Id-\Delta)^{s/2}\big(u(\tau,.) 
\otimes v(\tau,.)\big)\Big)_j =\notag \\
e^{(t - \tau)\Delta}\sum_{l, k = 1}^d\Big(\delta_{jk} - 
\frac{\partial_j\partial_k}{\Delta}\Big)\partial_l(Id-\Delta)^{s/2} 
\big(u_l(\tau,.)v_k(\tau,.)\big). \notag
\end{gather}
From the property of the Fourier transform we have
\begin{gather}
\Big(e^{(t - \tau)\Delta}\mathbb{P}\nabla .(Id-\Delta)^{s/2}
\big(u(\tau,.) \otimes v(\tau,.)\big)\Big)_j^\wedge(\xi)  = \notag \\  
e^{-(t - \tau)|\xi|^2}\sum_{l, k = 1}^d \Big(\delta_{jk} - 
\frac{\xi_j\xi_k}{|\xi|^2}\Big)(i\xi_l)\Big((Id-\Delta)^{s/2}
\big(u_l(\tau,.)v_k(\tau,.)\big)\Big)^\wedge(\xi), \notag
\end{gather}
and therefore
\begin{gather}
\Big(e^{(t - \tau)\Delta}\mathbb{P}\nabla .(Id-\Delta)^{s/2}
\big(u(\tau,.) \otimes v(\tau,.)\big)\Big)_j = \notag \\
 \frac{1}{(t - \tau)^{\frac{d + 1}{2}}}\sum_{l, k = 1}^d K_{l, k, j}
\Big(\frac{.}{\sqrt{t - \tau}}\Big)*\Big((Id-\Delta)^{s/2}
\big(u_l(\tau,.)v_k(\tau,.)\big)\Big), \label{2.1.3}
\end{gather}
where 
$$
\widehat{K_{l, k, j}}(\xi)= \frac{1}{(2\pi)^{d/2}}.e^{-|\xi|^2}
\Big(\delta_{jk} - \frac{\xi_j\xi_k}{|\xi|^2}\Big)(i\xi_l).
$$
Applying Proposition 11.1 (\cite{P. G. Lemarie-Rieusset 2002}, p. 107)  
with $|\alpha| = 1$ we obtain
$$
|K_{l, k, j}(x)| \lesssim \frac{1}{(1 + |x|)^{d + 1}}.
$$
Thus, the tensor $K(x) = \{K_{l, k, j}(x)\}$ satisfies 
\begin{equation}\label{2.1.4}
|K(x)| \lesssim \frac{1}{(1 + |x|)^{d + 1}}.
\end{equation}
So, we can rewrite the equality \eqref{2.1.3}  in the tensor form 
$$
e^{(t - \tau)\Delta}\mathbb{P}\nabla .(Id-\Delta)^{s/2}
\big(u(\tau,.) \otimes v(\tau,.)\big)  = 
$$
$$
\frac{1}{(t - \tau)^{\frac{d + 1}{2}}}
K\Big(\frac{.}{\sqrt{t - \tau}}\Big)*\Big((Id-\Delta)^{s/2}
\big(u(\tau,.) \otimes v(\tau,.)\big)\Big).
$$
Set
\begin{equation}\label{2.1.6}
\frac{1}{\tilde{p}} = \frac{2}{p} - \frac{s}{d},\ \frac{1}{h} 
= \frac{s}{d} - \frac{1}{p} + 1. 
\end{equation}
Note that from the inequalities \eqref{lem.2.1.2.a} and 
\eqref{lem.2.1.2.b}, we can check that the following relations are satisfied
$$
1 < h, \tilde{p} < \infty \ \text{and}\ \frac{1}{p} + 1 
= \frac{1}{h} + \frac{1}{\tilde{p}}.
$$
Applying the Young inequality for convolution we obtain
\begin{gather}
\Big\|e^{(t - \tau)\Delta}\mathbb{P}\nabla .(Id-\Delta)^{s/2}\big(u(\tau,.) 
\otimes v(\tau,.)\big)\Big\|_{L^p} \lesssim \notag \\
\frac{1}{(t - \tau)^{\frac{d + 1}{2}}}\Big\|K\Big(\frac{.}{\sqrt{t - \tau}}\Big)\Big\|_{L^h}
\Big\|(Id-\Delta)^{s/2}\big(u(\tau,.) \otimes v(\tau,.)\big)\Big\|_{L^{\tilde{p}}}.\label{2.1.5}
\end{gather}
Applying Lemma \ref{Them}
\begin{gather}
\Big\|(Id-\Delta)^{s/2}\big(u(\tau,.) \otimes v(\tau,.)\big)\Big\|_{L^{\tilde{p}}} 
= \big\|u(\tau,.)\otimes v(\tau,.)\big\|_{H^s_{{\tilde{p}}}} \notag \\
\lesssim \big\|u(\tau,.)\big\|_{H^s_p}\big\|v(\tau,.)\big\|_{H^s_{p}}.
\label{2.1.7}
\end{gather}
From the estimate \eqref{2.1.4} and the equality \eqref{2.1.6}, we have
\begin{gather}
\Big\|K\Big(\frac{.}{\sqrt{t - \tau}}\Big)\Big\|_{L^h} 
= (t - \tau)^{\frac{d}{2h}}\big\|K\big\|_{L^h} \simeq (t - \tau)^{\frac{s}{2} 
- \frac{d}{2p} + \frac{d}{2}}.
\label{2.1.8}
\end{gather}
The inequalities \eqref{2.1.5}, \eqref{2.1.7}, and \eqref{2.1.8} imply that 
\begin{gather}
\Big\|e^{(t - \tau)\Delta}\mathbb{P}\nabla .(Id-\Delta)^{s/2}\big(u(\tau,.) 
\otimes v(\tau,.)\big)\Big\|_{L^p} \lesssim \notag \\
 (t - \tau)^{\frac{s}{2} - \frac{d}{2p} 
- \frac{1}{2}}\big\|u(\tau,.)\big\|_{H^s_p}\big\|v(\tau,.)\big\|_{H^s_p}. \label{2.1.9}
\end{gather}
From the inequalities \eqref{2.1.2} and \eqref{2.1.9}, we get 
\begin{gather}
\big\|B(u, v)(t)\big\|_{H^s_p} \lesssim \int_0^t (t - \tau)^{\frac{s}{2} 
- \frac{d}{2p} - \frac{1}{2}}\big\|u(\tau,.)\big\|_{H^s_p}
\big\|v(\tau,.)\big\|_{H^s_p}\dif \tau  . \notag
\end{gather}
Applying  of Proposition 2.4 (c) in (\cite{P. G. Lemarie-Rieusset 2002}, p. 20) 
for the convolution in the Lorentz spaces, we have the following estimates
\begin{gather}
\Big\|\big\|B(u, v)(t)\big\|_{H^s_p}\Big\|_{L^r_t(0, T)}
= \Big\|\big\|B(u, v)(t)\big\|_{H^s_p}\Big\|_{L^{r, r}_t(0, T)} \notag \\
 \leq \Big\|\big\|B(u, v)(t)\big\|_{H^s_p}\Big\|_{L^{r, \frac{r}{2}}_t(0, T)} 
\label{2.1.10} \lesssim \notag \\ 
 \big\|1_{[0, T]}t^{\frac{s}{2} - \frac{d}{2p} - \frac{1}{2}}\big\|_{L^{r', \infty}} 
\Big\|\big\|u(t,.)\big\|_{H^s_p}\big\|v(t,.)\big\|_{H^s_p}
\Big\|_{L^{\frac{r}{2}, \frac{r}{2}}_t(0, T)},
\label{2.1.11}
\end{gather}
where $\frac{1}{r'} + \frac{1}{r} = 1$ and $1_{[0, T]}$ 
is the indicator function of set $[0, T]$ on $\mathbb R$. \\
By applying the H\"{o}lder inequality  we get 
\begin{gather}
\Big\|\big\|u(t,.)\big\|_{H^s_p}\big\|v(t,.)\big\|_{H^s_p}\Big\|_{L^{\frac{r}{2}, 
\frac{r}{2}}_t(0, T)} = 
\Big\|\big\|u(t,.)\big\|_{H^s_p}\big\|v(t,.)\big\|_{H^s_p}
\Big\|_{L^{\frac{r}{2}}_t(0, T)} \notag\\
\leq \Big\|\big\|u(t,.)\big\|_{H^s_p}\Big\|_{L^r_t(0, T)}
\Big\|\big\|v(t,.)\big\|_{H^s_p}\Big\|_{L^r_t(0, T)}.\label{2.1.12}
\end{gather}
Note that 
\begin{equation}\label{2.1.13}
\Big\|1_{[0, T]}t^{\frac{s}{2} - \frac{d}{2p} - \frac{1}{2}}
\Big\|_{L^{r', \infty}}
\simeq  T^{\frac{1}{2}(1 + s - \frac{2}{r} - \frac{d}{p})}.
\end{equation}
Therefore the inequality \eqref{2.1.1} can be deduced from 
the inequalities \eqref{2.1.11}, \eqref{2.1.12}, and \eqref{2.1.13}. \qed  
\begin{Not}
Lemma \ref{lem.2.1.2} is still valid 
when the inhomogeneous space $H^s_p$ 
is replaced by the homogeneous space $\dot{H}^s_p$.
\end{Not}
\begin{Lem}\label{lem.2.1.3}
Let 
$$
d \geq 3, \ 0 \leq s < d, \ p > 1, \ r > 2, and \ T>0
$$ 
be such that
$$
 \frac{1}{p} < \frac{1}{2} +\frac{s}{2d}, \frac{2}{p}\geq \frac{s+1}{d},
\ and\ \frac{2}{r} + \frac{d}{p} - s = 1.
$$
Then the bilinear operator $B(u, v)(t)$ is continuous from 
$$
L^r\big([0, T]; \dot{H}^s_p\big) 
\times L^r\big([0, T]; \dot{H}^s_p\big) 
$$
into
$$
L^\infty\Big([0, T]; \dot{B}^{\frac{d}{\tilde{p}} - 1, 
\frac{r}{2}}_{\tilde{p}}\Big),
$$
where
$$
\frac{1}{\tilde{p}} = \frac{2}{p} - \frac{s}{d},
$$
and we have the inequality
\begin{equation}\label{lem.2.1.3.1}
\big\|B(u, v)\big\|_{L^\infty\big([0, T]; \dot{B}^{\frac{d}{\tilde{p}} - 1, 
\frac{r}{2}}_{\tilde{p}}\big)} 
\leq C\big\|u\big\|_{L^r([0, T]; \dot{H}^s_p)}
\big\|v\big\|_{L^r([0, T]; \dot{H}^s_p)},
\end{equation}
where C is a positive constant independent of T.
\end{Lem}
\textbf{Proof}: To prove this lemma by duality (in the x-variable), 
$\big($see Proposition 3.9 in (\cite{P. G. Lemarie-Rieusset 2002}, p. 29)$\big)$, 
let us consider an arbitrary test function \linebreak $h(x) \in \mathcal{S}(\mathbb{R}^d)$  
and evaluate the quantity 
\begin{gather}
I_t =  \big<B(u,v)(t),h\big> = \int_{\mathbb{R}^d}
\big(B(u,v)(t)\big)(x)h(x)\dif x.\label{2.1.14}
\end{gather}
We have
\begin{gather}
\big<B(u,v)(t),h\big> = \int_0^t\big<e^{(t- \tau) \Delta} \mathbb{P} 
\nabla .\big(u(\tau,.)\otimes v(\tau,.)\big),h\big>\dif \tau  =\notag \\
\int_0^t\Big<e^{(t- \tau) \Delta}\dot{\Lambda}\mathbb{P} 
\frac{\nabla}{\dot{\Lambda}}.\big(u(\tau,.)\otimes v(\tau,.)\big),h\Big>\dif \tau   = \notag \\
\int_0^t\Big<\mathbb{P} \frac{\nabla}{\dot{\Lambda}} .\big(u(\tau,.)
\otimes v(\tau,.)\big) ,e^{(t- \tau) \Delta}\dot{\Lambda}h\Big>\dif \tau  = \notag \\
\int_0^t\Big<\mathbb{P} \frac{\nabla}{\dot{\Lambda}}.
\dot{\Lambda}^s\big(u(\tau,.)\otimes v(\tau,.)\big) ,e^{(t- \tau) \Delta}
\dot{\Lambda}\dot{\Lambda}^{-s}h\Big>\dif \tau  . \label{2.1.15}
\end{gather}
By applying the H\"{o}lder inequality in the x-variable, 
from the equality \eqref{2.1.15} and the fact that (see \cite{P. G. Lemarie-Rieusset 2002})
$$
\mathbb{P}\ \text{and}\ \frac{\nabla}{\dot{\Lambda}} 
\ \text{are continuous from} \ L^p \ {\rm into} 
\ L^p, 1 < p <\infty, 
$$
we get
\begin{gather}
|I_t| \leq \int_0^t\Big\|\mathbb{P} \frac{\nabla}{\dot{\Lambda}}. 
\dot{\Lambda}^s\big(u(\tau,.)
\otimes v(\tau,.)\big)\Big\|_{L^{\tilde{p}}}\big\|e^{(t- \tau) \Delta}
\dot{\Lambda}\dot{\Lambda}^{-s}h\big\|_{L^{\tilde{p}'}}\ \dif \tau  \notag \\
\lesssim \int_0^t\big\|\dot{\Lambda}^s\big(u(\tau,.)
\otimes v(\tau,.)\big)\big\|_{L^{\tilde{p}}}\big\|e^{(t- \tau) \Delta}
\dot{\Lambda}\dot{\Lambda}^{-s}h\big\|_{L^{\tilde{p}'}}\ \dif \tau , \label{2.1.16}
\end{gather}
where 
$$
\frac{1}{\tilde{p}} + \frac{1}{\tilde{p}'} =1.
$$
Applying Lemma \ref{Them}, we have
\begin{gather}
\big\|\dot{\Lambda}^s\big(u(\tau,.)\otimes v(\tau,.)\big)\big\|_{L^{\tilde{p}}} 
=\big\|u(\tau,.)\otimes v(\tau,.)\big\|_{\dot{H}_{\tilde{p}}^s}\notag\\
 \lesssim  \big\|u(\tau,.)\big\|_{\dot{H}^s_p} \big\|v(\tau,.)
\big\|_{\dot{H}^s_p}.\label{2.1.17}
\end{gather}
From the inequalities \eqref{2.1.16} and \eqref{2.1.17}, 
applying the H\"{o}lder inequality in the t-variable, we deduce that 
\begin{gather}
|I_t| \lesssim \int_0^t \big\|u(\tau,.)\big\|_{\dot{H}^s_p} 
\big\|v(\tau,.)\big\|_{\dot{H}^s_p}\big\|e^{(t- \tau) \Delta}
\dot{\Lambda}\dot{\Lambda}^{-s}h\big\|_{L^{{\tilde{p}'}}}\ \dif \tau  \leq \notag \\
\Big(\int_0^t \big(\big\|u(\tau,.)\big\|_{\dot{H}^s_p} 
\big\|v(\tau,.)\big\|_{\dot{H}^s_p}\big)^{\frac{r}{2}}\dif \tau \Big)^{\frac{2}{r}}
\Big(\int_0^t\big(\big\|e^{(t- \tau) \Delta}\dot{\Lambda}\dot{\Lambda}^{-s}
h\big\|_{L^{{\tilde{p}'}}}\big)^{\frac{r}{r-2}}d \tau \Big)^{\frac{r-2}{r}} \notag \\
\leq \big\|u\big\|_{L^r([0, T]; \dot{H}^s_p)}
\big\|v\big\|_{L^r([0, T]; \dot{H}^s_p)}\Big(\int_0^t\big(\big\|e^{(t- \tau) \Delta}
\dot{\Lambda}\dot{\Lambda}^{-s}h\big\|_{L^{{\tilde{p}'}}}\big)^{\frac{r}{r-2}}
\dif \tau \Big)^{\frac{r-2}{r}}. \label{2.1.18}
\end{gather}
From Lemma \ref{lem.2.1.4} and note that $\dot{\Lambda}^{s_0}$ 
is an isomorphism from $\dot{B}^{s,p}_q$ to $\dot{B}^{s-s_0,p}_q$ 
(see \cite{Bjorn Jawerth: 1977}),  we have the following estimates
\begin{gather}
\Big(\int_0^t\big(\big\|e^{(t- \tau) \Delta}\dot{\Lambda}
\dot{\Lambda}^{-s}h\big\|_{L^{{\tilde{p}'}}}\big)^{\frac{r}{r-2}}
\dif \tau \Big)^{\frac{r-2}{r}} \leq 
 \Big(\int_0^\infty\big(\big\|e^{t \Delta}\dot{\Lambda}
\dot{\Lambda}^{-s}h\big\|_{L^{{\tilde{p}'}}}\big)^{\frac{r}{r-2}}
\dif t\big)^{\frac{r-2}{r}} \notag \\
= \Big(\int_0^\infty\big(t^{\frac{r-4}{2r}}
\big\|e^{t \Delta}t^{\frac{1}{2}}\dot{\Lambda}\dot{\Lambda}^{-s}
h\big\|_{L^{{\tilde{p}'}}}\big)^{\frac{r}{r-2}}\frac{\dif t}{t}\Big)^{\frac{r-2}{r}}  
\simeq \big\|\dot{\Lambda}^{-s}h\big\|_{\dot{B}_{\tilde{p}'}^{\frac{4-r}{r}, 
{\frac{r}{r-2}}}} \notag\\
 \simeq  \big\|h\big\|_{\dot{B}_{\tilde{p}'}^{\frac{4-r}{r}-s, {\frac{r}{r-2}}}} 
= \big\|h\big\|_{\dot{B}_{\tilde{p}'}^{1 - \frac{d}{\tilde{p}}, 
{\frac{r}{r-2}}}}  . \label{2.1.19}
\end{gather}
From the equality \eqref{2.1.14} and the inequalities  \eqref{2.1.18} 
and \eqref{2.1.19}, we get 
\begin{gather}
\big|\big<B(u,v)(t),h\big>\big| \lesssim  \big\|u\big\|_{L^r([0, T]; \dot{H}^s_p)}
\big\|v\big\|_{L^r([0, T]; \dot{H}^s_p)}
\Big\|h\Big\|_{\dot{B}_{\tilde{p}'}^{1 - \frac{d}{\tilde{p}}, 
{\frac{r}{r-2}}}}.\notag
\end{gather}
However, $\dot{B}_{\tilde{p}'}^{1 - \frac{d}{\tilde{p}}, {\frac{r}{r-2}}}$ 
is exactly the dual of $\dot{B}_{\tilde{p}}^{\frac{d}{\tilde{p}} - 1, \frac{r}{2}}$,  
(the restriction $\frac{2}{p}\geq\frac{s+1}{d}$ is mainly because we are interested 
in non-negative indexes), 
therefore we conclude that 
\begin{equation}\label{lem.2.1.3.2}
\Big\|B(u,v)(t)\Big\|_{\dot{B}_{\tilde{p}}^{\frac{d}{\tilde{p}} - 1, \frac{r}{2}}} \lesssim 
 \big\|u\big\|_{L^r([0, T]; \dot{H}^s_p)}
\big\|v\big\|_{L^r([0, T]; \dot{H}^s_p)}, \ 0 \leq t \leq T.
\end{equation}
Finally, the estimate \eqref{lem.2.1.3.1} can be deduced from 
the inequality \eqref{lem.2.1.3.2}. \qed \\

Combining Theorem \ref{Th.2.1} with Lemma \ref{lem.2.1.2}, 
we get the following existence results,  the particular case of which, when $s=0$, 
was obtained in \cite{P. G. Lemarie-Rieusset 2002}.
\begin{Th}\label{th.3.3.1}
Let 
$$
d \geq 3, s \geq 0, \ p > 1, \ and\ r > 2,
$$
be such that
$$
\frac{s}{d} < \frac{1}{p} < \frac{1}{2} +\frac{s}{2d}
\ and\  \frac{2}{r} + \frac{d}{p} - s  \leq 1.
$$
{\rm (a)} There exists a positive constant $\delta_{s,p,r,d}$ 
such that for all $T > 0$ and for all $u_0 \in \mathcal{S}'(\mathbb{R}^d)$ 
with ${\rm div}(u) = 0$, satisfying
\begin{equation}\label{3.3.1.2.a}
T^{\frac{1}{2}(1 + s - \frac{2}{r} -  \frac{d}{p})}
\big\|e^{t\Delta}u_0\big\|_{L^r([0, T]; \dot{H}^s_p)} \leq \delta_{s,p,r,d},
\end{equation}
there is a unique mild solution $u \in L^r\big([0, T]; \dot{H}^s_p)\big)$ for NSE. \\
If 
$$
e^{t\Delta}u_0 \in L^r\big([0, 1]; \dot{H}^s_p\big),
$$
 then the inequality \eqref{3.3.1.2.a} holds when $T(u_0)$ is small enough.\\ 
{\rm (b)} If \ $\frac{2}{r} + \frac{d}{p} - s = 1$ then there exists  
a positive constant $\delta_{s,p,d}$ such that we can take  $T = \infty$ whenever $\big\|e^{t\Delta}u_0\big\|_{L^r([0, \infty]; \dot{H}^s_p)} \leq \delta_{s,p,d}$.
\end{Th}
\textbf{Proof}: (a) From Lemma \ref{lem.2.1.2}, we use the estimate 
$$
\big\|B\big\|_{L^r([0, T]; \dot{H}^s_p)} 
\leq C_{s,p,r,d}T^{\frac{1}{2}(1 + s - \frac{2}{r} -  \frac{d}{p})},
$$ 
where $C_{s,p,r,d}$ is a positive constant independent of $T$. 
From Theorem \ref{Th.2.1} and the above inequality, 
we deduce  the existence of a solution to the 
Navier-Stokes equations on the interval $(0, T)$ with  
$$
4C_{s,p,r,d}T^{\frac{1}{2}(1 + s  - \frac{2}{r} 
- \frac{d}{p})}\big\|e^{t\Delta}u_0\big\|_{L^r([0, T]; \dot{H}^s_p)} \leq 1.
$$
If $e^{t\Delta}u_0 \in L^r([0, 1]; \dot{H}^s_p)$ then this condition 
is fulfilled for $T = T(u_0)$ small enough, this is obvious for the 
case when $\frac{2}{r} + \frac{d}{p} - s < 1$ 
since \linebreak $\underset{T \rightarrow 0}{\rm lim}
T^{\frac{1}{2}(1 + s - \frac{2}{r} - \frac{d}{p})} = 0$. 
For the case when $\frac{2}{r} + \frac{d}{p} - s  =1$, 
the condition is fulfilled since we have $\underset{T \rightarrow 0}{\rm lim}\big\|e^{t\Delta}u_0\big\|_{L^r([0, T]; \dot{H}^s_p)} = 0$.\\
(b) This is obvious. \qed 
\begin{Not}
From Theorem 5.3 (\cite{P. G. Lemarie-Rieusset 2002}, p. 44), 
if $u_0 \in B_p^{s - \frac{2}{r},r}\cap \mathcal{S}'(\mathbb{R}^d)$ then  
$e^{t\Delta}u_0 \in L^r([0, 1]; \dot{H}^s_p)$. 
From Lemma \ref{lem.2.1.4'}, if $u_0\in \mathcal{S}'(\mathbb{R}^d)$ 
the two quantities 
$\big\|u_0\big\|_{\dot{B}_p^{s -\frac{2}{r},r}}$ 
and $\big\|e^{t\Delta}u_0\big\|_{L^{r}([0, \infty]; \dot{H}^s_p)}$
 are equivalent.
\end{Not}

\vskip 0,3cm
\textbf{3.2. Solutions to the Navier-Stokes equations with initial 
value in the critical spaces} $H^{\frac{d}{q} - 1}_q(\mathbb{R}^d)$ 
\textbf{and} $\dot{H}^{\frac{d}{q} - 1}_q(\mathbb{R}^d)$ \textbf{for} 
$3 \leq d \leq 4, \ 2 \leq q \leq d$.
\vskip 0.4cm

\begin{Lem}\label{lem.2.2.2}
Let  $d \geq 3\ and \  2 \leq  q \leq d$. Then the bilinear 
operator $B(u, v)(t)$ is continuous from 
$$
L^4\Big([0, T]; \dot{H}^{\frac{d}{q} - 1}_{\frac{2dq}{2d - q}}\Big)
 \times L^4\Big([0, T]; 
\dot{H}^{\frac{d}{q} - 1}_{\frac{2dq}{2d - q}}\Big)
$$
into
$$
 L^\infty\Big([0, T]; \dot{B}^{\frac{d}{q} - 1,2}_q\Big),
$$
and we have the inequality
\begin{gather}
\big\|B(u, v)\big\|_{L^\infty([0, T]; \dot{H}^{\frac{d}{q} - 1}_q)} 
\lesssim \big\|B(u, v)\big\|_{L^\infty\big([0, T]; 
\dot{B}^{\frac{d}{q} - 1,2}_q\big)} \notag \\ 
\leq  C\big\|u\big\|_{L^4\big([0, T]; 
\dot{H}^{\frac{d}{q} - 1}_{\frac{2dq}{2d - q}}\big)}
\big\|v\big\|_{L^4\big([0, T]; 
\dot{H}^{\frac{d}{q} - 1}_{\frac{2dq}{2d - q}}\big)}, \label{2.2.1}
\end{gather}
where C is a positive constant and independent of T.
\end{Lem}
\textbf{Proof}: Applying Lemma \ref{lem.2.1.3} 
with $r = 4, p = \frac{2dq}{2d - q},\ {\rm and}\  s = \frac{d}{q} - 1$, we get
\begin{gather}
\frac{1}{\tilde{p}} = \frac{2}{p} - \frac{s}{d} 
= \frac{2d - q}{dq} - \frac{\frac{d}{q} - 1}{d} =\frac{1}{q}, \notag\\
\big\|B(u, v)\big\|_{L^\infty\big([0, T]; \dot{B}^{\frac{d}{q} - 1, 2}_{q}\big)} 
\lesssim \big\|u\big\|_{L^4\big([0, T]; 
\dot{H}^{\frac{d}{q} - 1}_{\frac{2dq}{2d - q}}\big)}
\big\|v\big\|_{L^4\big([0, T]; 
\dot{H}^{\frac{d}{q} - 1}_{\frac{2dq}{2d - q}}\big)}.\label{2.2.2}
\end{gather}
From (b) of Lemma \ref{lem.2.2.3}, we have 
\begin{equation}\label{2.2.4}
\dot{B}^{\frac{d}{q} -1, 2}_q \hookrightarrow  \dot{H}^{\frac{d}{q} - 1}_q.
\end{equation}
Finally, the estimate \eqref{2.2.1} can be deduced 
from the inequality \eqref{2.2.2} and the imbedding \eqref{2.2.4}. \qed 
\begin{Lem}\label{lem.2.2.5}
Let  $d \geq 3\ and \  2 \leq  q \leq d$. 
Then the bilinear operator $B(u, v)(t)$ is continuous from 
$$
L^4\Big([0, T]; H^{\frac{d}{q} - 1}_{\frac{2dq}{2d - q}}\Big)
 \times L^4\Big([0, T]; 
H^{\frac{d}{q} - 1}_{\frac{2dq}{2d - q}}\Big)
$$
into
$$
L^\infty\Big([0, T]; H^{\frac{d}{q} - 1}_q\Big),
$$
and we have the inequality
\begin{gather}
\big\|B(u, v)\big\|_{L^\infty\big([0, T]; H^{\frac{d}{q} - 1}_q\big)} 
\leq  C\big\|u\big\|_{L^4\big([0, T]; H^{\frac{d}{q} - 1}_{\frac{2dq}{2d - q}}\big)}
\big\|v\big\|_{L^4\big([0, T]; 
H^{\frac{d}{q} - 1}_{\frac{2dq}{2d - q}}\big)}, \label{2.2.7}
\end{gather}
where C is a positive constant and independent of T.
\end{Lem}
\textbf{Proof}: To prove this lemma by duality (in the x-variable), 
let us consider an arbitrary test function $h(x) \in \mathcal{S}(\mathbb{R}^d)$. 
Similar to the proof of Lemma \ref{lem.2.1.3}, we have 
\begin{gather}
\Big|\big<(\sqrt{Id-\Delta})^{\frac{d}{q} - 1}B(u,v)(t),h\big>\Big| \lesssim 
 \big\|u\big\|_{L^4\big([0, T]; H^{\frac{d}{q} - 1}_{\frac{2dq}{2d - q}}\big)}
\big\|v\big\|_{L^4\big([0, T]; H^{\frac{d}{q} - 1}_{\frac{2dq}{2d - q}}
\big)}\big\|h\big\|_{\dot{B}_{q'}^{0, 2}},\notag
\end{gather}
where
$$
\frac{1}{q} + \frac{1}{q'} = 1.
$$
However the dual space of $\dot{B}_{q'}^{0, 2}$ is $\dot{B}_{q}^{0, 2}$,
therefore we get
\begin{gather}
\Big\|(\sqrt{Id-\Delta})^{\frac{d}{q} - 1}B(u,v)(t)\Big\|_{\dot{B}_{q}^{0, 2}} \lesssim 
 \big\|u\big\|_{L^4\big([0, T]; H^{\frac{d}{q} - 1}_{\frac{2dq}{2d - q}}
\big)}\big\|v\big\|_{L^4\big([0, T]; 
H^{\frac{d}{q} - 1}_{\frac{2dq}{2d - q}}\big)}. \label{2.2.8}
\end{gather}
 From (b) of Lemma \ref{lem.2.2.3} and the estimate \eqref{2.2.8}, we have
\begin{gather}
\Big\|B(u,v)(t)\Big\|_{H^{\frac{d}{q} - 1}_q}
=\Big\|(\sqrt{Id-\Delta})^{\frac{d}{q} - 1}B(u,v)(t)\Big\|_{L^q} = \notag \\
 \Big\|(\sqrt{Id-\Delta})^{\frac{d}{q} - 1}B(u,v)(t)\Big\|_{\dot{H}^0_q} 
\lesssim  \Big\|(\sqrt{Id-\Delta})^{\frac{d}{q} - 1}
B(u,v)(t)\Big\|_{\dot{B}^{0, 2}_q} \notag \\
 \lesssim \big\|u\big\|_{L^4\big([0, T]; 
H^{\frac{d}{q} - 1}_{\frac{2dq}{2d - q}}\big)}
\big\|v\big\|_{L^4\big([0, T]; H^{\frac{d}{q} - 1}_{\frac{2dq}{2d - q}}\big)}
\label{2.2.9}, \ 0 \leq t \leq T.
\end{gather}
Finally, the estimate \eqref{2.2.7} can be deduced 
from the inequality \eqref{2.2.9}. \qed 

\begin{Lem}\label{lem.2.2.4}Let $d \geq 3$ and $2 \leq q \leq 4$.\\
{\rm (a)} If \ $u_0 \in H^{\frac{d}{q} - 1}_q(\mathbb{R}^d)$ then 
$$
\big\|e^{t\Delta}u_0\big\|_{L^4\big([0, \infty); H^{d/q - 1}_{2dq/(2d - q)}\big)} 
\lesssim \big\|u_0\big\|_{H^{d/q - 1}_q}.
$$
{\rm (b)} If \ $u_0 \in \dot{H}^{\frac{d}{q} - 1}_q(\mathbb{R}^d)$ then 
$$
\Big\|e^{t\Delta}u_0\Big\|_{L^4\big([0, \infty); \dot{H}^{d/q - 1}_{2dq/(2d - q)}\big)} 
\simeq  \big\|u_0\big\|_{\dot{B}_{2dq/(2d - q)}^{d/q -3/2, 4}} 
\lesssim \big\|u_0\big\|_{\dot{H}^{d/q - 1}_q}.
$$
\end{Lem}
\textbf{Proof}: (a) From Lemma \ref{lem.2.1.4}, we have the estimates 
\begin{gather}
\Big\|e^{t\Delta}u_0\Big\|_{L^4\big([0, \infty); H^{d/q - 1}_{2dq/(2d - q)}\big)} 
=\Big(\int_0^\infty \Big\|e^{t\Delta}(\sqrt{Id -\Delta})^{d/q - 1}
u_0\Big\|^4_{L^{2dq/(2d - q)}}\dif t\Big)^{1/4}  \notag \\
=\Big(\int_0^\infty \Big(t^{\frac{1}{4}}\Big\|e^{t\Delta}
(\sqrt{Id -\Delta})^{d/q - 1}u_0\Big\|_{L^{2dq/(2d - q)}}\Big)^4
\frac{\dif t}{t}\Big)^{1/4} \notag \\
\simeq \Big\|(\sqrt{Id -\Delta})^{d/q - 1}
u_0\Big\|_{\dot{B}_{2dq/(2d - q)}^{-1/2, 4}}. \label{2.2.5}
\end{gather}
Applying (b), (c), and (d) of Lemma \ref{lem.2.2.3} in order to obtain
\begin{equation}\label{2.2.6}
L^q = \dot{H}_q^0  \hookrightarrow \dot{B}_q^{0, {q}}  \hookrightarrow
 \dot{B}_q^{0, 4} \hookrightarrow \dot{B}_{2dq/(2d - q)}^{-1/2, 4}.
\end{equation}
From the inequality \eqref{2.2.5} and the imbedding \eqref{2.2.6}, we get 
\begin{gather*}
\Big\|e^{t\Delta}u_0\Big\|_{L^4\big([0, \infty); H^{d/q - 1}_{2dq/(2d - q)}\big)} 
\simeq \Big\|(\sqrt{Id -\Delta})^{d/q - 1}u_0\Big\|_{\dot{B}_{2dq/(2d - q)}^{-1/2, 4}}\\
 \lesssim  \big\|(\sqrt{Id -\Delta})^{d/q - 1}u_0\big\|_{L^q}
= \big\|u_0\big\|_{H_{q}^{d/q - 1}}.
\end{gather*}
(b) Similar to the proof of (a) we have
\begin{gather}
\big\|e^{t\Delta}u_0\big\|_{L^4\big([0, \infty); 
\dot{H}^{d/q - 1}_{2dq/(2d - q)}\big)} 
\simeq \big\|\dot{\Lambda}^{\frac{d}{q} - 1}
u_0\big\|_{\dot{B}_{2dq/(2d - q)}^{-1/2, 4}} \notag \\
\lesssim \big\|\dot{\Lambda}^{\frac{d}{q} - 1}u_0\big\|_{L^q} 
= \big\|u_0\big\|_{\dot{H}_{q}^{d/q - 1}},\notag \\
\text{and} \ \ \  \Big\|\dot{\Lambda}^{\frac{d}{q} - 1}
u_0\Big\|_{\dot{B}_{2dq/(2d - q)}^{-1/2, 4}} 
\simeq  \big\|u_0\big\|_{\dot{B}_{2dq/(2d - q)}^{d/q -3/2, 4}} . \ \ \ \ \qed \notag
\end{gather}
Combining Theorem \ref{Th.2.1} with Lemmas \ref{lem.2.2.6},
  \ref{lem.2.1.2}, \ref{lem.2.2.2}, and \ref{lem.2.2.4} 
we obtain the following existence result.
\begin{Th}\label{Th.1}
Let $3 \leq d \leq 4\ and \ 2 \leq q \leq d$. There exists a positive constant $\delta_{q, d}$ 
such that for all $T > 0$ and for all $u_0 \in \dot{H}^{d/q -1}_q(\mathbb{R}^d)$ 
with ${\rm div}(u_0) = 0$ satisfying
\begin{equation}\label{2.2.10}
\big\|e^{t\Delta}u_0\big\|_{L^4\big([0, T]; 
\dot{H}^{d/q - 1}_{2dq/(2d - q)}\big)} \leq \delta_{q, d},
\end{equation}
NSE has a unique mild solution $u \in L^4\big([0, T]; 
\dot{H}^{d/q - 1}_{2dq/(2d - q)}\big) \cap C\big([0, T]; 
\dot{H}^{d/q - 1}_q\big)$. Denoting $w = u - e^{t\Delta}u_0$, then we have
$$
w \in L^4\big([0, T]; 
\dot{H}^{d/q - 1}_{2dq/(2d - q)}\big) \cap L^\infty\big([0, T]; 
\dot{B}^{d/q - 1, 2}_q\big).
$$
Finally, we have 
\begin{gather*}
\big\|e^{t\Delta}u_0\big\|_{L^4\big([0, T]; \dot{H}^{d/q - 1}_{2dq/(2d - q)}\big)} 
\lesssim \big\|u_0\big\|_{\dot{B}^{d/q - 3/2, 4}_{2dq/(2d - q)}}
 \lesssim \big\|u_0\big\|_{\dot{H}^{d/q - 1}_q}, 
\end{gather*}
in particular,  for arbitrary 
$u_0 \in \dot{H}^{d/q -1}_q(\mathbb{R}^d)$ the inequality \eqref{2.2.10} 
holds when $T(u_0)$ is small enough; 
and there exists a positive constant $\sigma_{q, d}$ 
such that for all $\Big\|u_0\Big\|_{\dot{B}^{d/q - 3/2, 4}_{2dq/(2d - q)}} \leq \sigma_{q, d}$
we can take  $T = \infty$.
\end{Th}
\textbf{Proof}: By applying Lemma \ref{lem.2.1.2} 
with $r = 4, \ p = \frac{2dq}{2d - q},\ s = \frac{d}{q} - 1$, 
and notice that $1 + s - \frac{2}{r} - \frac{d}{p} = 0$ we have  
$$
\big\|B\big\|_{L^4\big([0, T]; 
\dot{H}^{d/q - 1}_{2dq/(2d - q)}\big)} \leq C_{q,d},
$$ 
where $C_{q,d}$ is a positive constant independent of  $T$. From 
Theorem \ref{Th.2.1} and the above inequality, we 
deduce that for any $u_0 \in \dot{H}^{\frac{d}{q} -1}_q$ such that
$$
{\rm div}(u_0)=0,\ \ \big\|e^{t\Delta}u_0\big\|_{L^4\big([0, T]; 
\dot{H}^{d/q - 1}_{2dq/(2d - q)}\big)} \leq \frac{1}{4C_{q,d}},
$$ 
NSE has a mild solution $u$ on the interval $(0, T)$ so that 
\begin{equation}\label{2.2.11}
u \in L^4\big([0, T]; \dot{H}^{d/q - 1}_{2dq/(2d - q)}\big).
\end{equation}
From Lemma \ref{lem.2.2.2} and \eqref{2.2.11}, we have 
$
B(u,u) \in L^\infty\big([0, T]; \dot{H}^{d/q - 1}_q\big).  
$
From (2) of Lemma \ref{lem.2.2.6}, we have 
$e^{t\Delta}u_0 \in L^\infty\big([0, T]; \dot{H}^{d/q - 1}_q\big)$. 
Therefore
$$
u = e^{t\Delta}u_0 - B(u,u) \in  L^\infty([0, T]; \dot{H}^{d/q - 1}_q).
$$
In the space $H^{d/2-1}$ or $L^d$ (see \cite{P. G. Lemarie-Rieusset 2002}), 
the solutions can also be constructed by a successive approximation via 
the integral equation and therefore they are continuous in time up to the 
initial time. Since $e^{t\Delta}$ is a $(C_0)$-semigroup 
in $H^s_q$ and $\dot H^s_q$ with finite integral-exponent $(q < \infty)$, 
by the same way as, we can easily show that the obtained mild solution $u \in C\big([0, T]; \dot{H}^{d/q - 1}_q\big)$. \\
From (b) of Lemma \ref{lem.2.2.4}, we have
\begin{gather*}
\big\|e^{t\Delta}u_0\big\|_{L^4\big([0, T]; \dot{H}^{d/q - 1}_{2dq/(2d - q)}\big)} 
\lesssim\big\|e^{t\Delta}u_0\big\|_{L^4\big([0, \infty); 
\dot{H}^{d/q - 1}_{2dq/(2d - q)}\big)}\notag \\
\simeq \Big\|u_0\Big\|_{\dot{B}^{d/q - 3/2, 4}_{2dq/(2d - q)}}
 \lesssim \big\|u_0\big\|_{\dot{H}^{d/q - 1}_q} <  \infty.
\end{gather*}
Hence, the left-hand side of the inequality \eqref{2.2.10} converges to $0$ when $T$ 
tends to $0$. Therefore, for arbitrary $u_0 \in \dot{H}^{\frac{d}{q} -1}_q$ 
there is $T(u_0)$ small enough  such that the inequality \eqref{2.2.10} holds. 
Also, there exists a positive constants $\sigma_{q, d}$ 
such that  for all 
$\Big\|u_0\Big\|_{\dot{B}^{d/q - 3/2, 4}_{2dq/(2d - q)}} \leq  \sigma_{q, d}$
and $T = \infty$ the inequality \eqref{2.2.10} holds.\qed 
\begin{Not} 
Theorem \ref{Th.1} in the particular case $q = d$ is  
Proposition 20.1 \linebreak in \cite{P. G. Lemarie-Rieusset 2002}. 
\end{Not}
\begin{Th}\label{Th.4}
Let $3 \leq d \leq 4\ and\ 2 \leq q \leq d$. There exists a positive 
constant $\delta_{q, d}$ such that for all $T > 0$ 
and for all $u_0 \in H^{\frac{d}{q} -1}_q(\mathbb{R}^d)$ 
with ${\rm div}(u_0) = 0$ satisfying  
\begin{equation}\label{2.2.15}
\big\|e^{t\Delta}u_0\big\|_{L^4\big([0, T]; H^{d/q - 1}_{2dq/(2d - q)}\big)} 
\leq \delta_{q, d},
\end{equation}
NSE has a unique mild solution  
$u \in L^4\big([0, T]; H^{d/q - 1}_{2dq/(2d - q)}\big) 
\cap C\big([0, T]; H^{d/q - 1}_q\big)$. Finally, we have 
$$
\big\|e^{t\Delta}u_0\big\|_{L^4\big([0, T]; H^{d/q - 1}_{2dq/(2d - q)}\big)} 
\leq  \big\|u_0\big\|_{H^{d/q -1}_q},
$$
in particular, for arbitrary $u_0 \in H^{\frac{d}{q} -1}_q$ 
the inequality \eqref{2.2.15} holds when $T(u_0)$ is small enough; 
\end{Th}
\textbf{Proof}:  The proof of Theorem \ref{Th.4} is 
similar to the one of Theorem  \ref{Th.1}, 
by combining Theorem \ref{Th.2.1} with Lemmas 
\ref{lem.2.2.6}, \ref{lem.2.1.2}, \ref{lem.2.2.5}, 
and  \ref{lem.2.2.4}. \qed 
\vskip 0.5cm
\textbf{3.3. Solutions to the Navier-Stokes equations 
with initial value in the critical spaces} 
$\dot{H}^{\frac{d}{q} - 1}_q(\mathbb{R}^d)$ \textbf{for} $d \geq 3
\ {\rm \textbf{and}}\ 1 < q \leq d$. 
\vskip 0.2cm
We consider two cases $2 < q \leq d$ and $1 < q \leq 2$ separately.
\vskip 0.2cm
\textbf{3.3.1. Solutions to the Navier-Stokes equations 
with initial value in the critical spaces} 
$\dot{H}^{\frac{d}{q} - 1}_q(\mathbb{R}^d)$ \textbf{for} $d \geq 3
\ {\rm \textbf{and}}\ 2 < q \leq d$. 
\vskip 0.2cm
\begin{Lem}\label{lem.2.3.2}
Let  $d \geq 3\ and \ 2 < q \leq  d $. Then for all $p$ such that 
$$
2 < p < {\rm min}\Big\{\frac{(d-2) q}{d-q}, \ d+ 2\Big\}, 
(if\ q  = d\ then\ \frac{(d-2) q}{d-q} = +\infty),
$$
the bilinear operator $B(u, v)(t)$ is continuous from 
$$
L^p([0, T]; \dot{H}^{\frac{2+d-p}{p}}_p)
\times L^p([0, T]; \dot{H}^{\frac{2+d-p}{p}}_p)
$$
into
$$
L^\infty\big([0, T]; \dot{B}^{\frac{d+p-2}{p}-1, 
\frac{p}{2}}_{\frac{d p}{d+p-2}}\big),
$$
and we have the inequality
\begin{gather}
\big\|B(u, v)\big\|_{L^\infty\big([0, T]; \dot{H}^{\frac{d}{q} - 1}_{q}\big)} \lesssim 
\big\|B(u, v)\big\|_{L^\infty\big([0, T]; \dot{B}^{\frac{d+p-2}{p}-1, 
\frac{p}{2}}_{\frac{d p}{d+p-2}}\big)} \notag \\
\leq C\big\|u\big\|_{L^p\big([0, T]; \dot{H}^{\frac{2+d-p}{p}}_p\big)}
\big\|v\big\|_{L^p\big([0, T]; \dot{H}^{\frac{2+d-p}{p}}_p\big)} \label{lem.2.3.2.1},
\end{gather}
where C is a positive constant independent of T.
\end{Lem}

\textbf{Proof}: Applying Lemma \ref{lem.2.1.3} 
with $r = p\ {\rm and}\ s = \frac{2+d-p}{p}$, we get
\begin{gather}
\frac{1}{\tilde{p}} = \frac{2}{p} - \frac{s}{d}
 =  \frac{d+p -2}{d p}, \notag\\
\Big\|B(u, v)\Big\|_{L^\infty\big([0, T]; \dot{B}^{\frac{d+p-2}{p}-1, 
\frac{p}{2}}_{\frac{d p}{d+p-2}}\big)} 
\lesssim \big\|u\big\|_{L^p\big([0, T]; \dot{H}^{\frac{2+d-p}{p}}_p\big)}
\big\|v\big\|_{L^p\big([0, T]; \dot{H}^{\frac{2+d-p}{p}}_p\big)}.\label{lem.2.3.2.2}
\end{gather}
Applying (e), (d), and (h) of Lemma \ref{lem.2.2.3} in order to obtain
\begin{equation}\label{lem.2.3.2.3}
\dot{B}^{\frac{d+p-2}{p}-1, \frac{p}{2}}_{\frac{d p}{d+p-2}} 
\hookrightarrow\dot{F}^{\frac{d+p-2}{p}-1, \frac{p}{2}}_{\frac{d p}{d+p-2}} 
 \hookrightarrow  \dot{F}_q^{\frac{d}{q}-1, 2}  =  \dot{H}_q^{\frac{d}{q}-1}.    
\end{equation}
Therefore the estimate \eqref{lem.2.3.2.1} is deduced 
from the inequality \eqref{lem.2.3.2.2} and the imbedding \eqref{lem.2.3.2.3}.
\begin{Lem}\label{lem.2.3.3}
Let  $2 <  q < p < + \infty$. Then for all $u_0 \in \dot{H}_q^{\frac{d}{q} - 1}$  
we have the estimates
$$
\big\|e^{t\Delta}u_0\big\|_{L^p\big([0, \infty); \dot{H}^{\frac{2+d-p}{p}}_p\big)} 
\simeq \big\|u_0\big\|_{\dot{B}_p^{\frac{d}{p}-1, p}}
\lesssim \big\|u_0\big\|_{\dot{H}_q^{\frac{d}{q} - 1}}.
$$
\end{Lem}
\textbf{Proof}:  From Lemma \ref{lem.2.1.4}, we have the estimates 
\begin{gather}
\big\|e^{t\Delta}u_0\big\|_{L^p\big([0, \infty); \dot{H}^{\frac{2+d-p}{p}}_p\big)}  
\simeq \big\|u_0\big\|_{\dot{B}_p^{\frac{d}{p}-1, p}} \ .\label{2.3.1}
\end{gather}
Applying (b), (d),  and (c) of Lemma \ref{lem.2.2.3} in order to obtain 
\begin{equation}\label{2.3.2}
\dot{H}_q^{\frac{d}{q} - 1} \hookrightarrow  \dot{B}_q^{\frac{d}{q} - 1,q}
 \hookrightarrow  \dot{B}_p^{\frac{d}{p} - 1,q}
 \hookrightarrow  \dot{B}_p^{\frac{d}{p} - 1,p}. 
\end{equation}
From the estimate \eqref{2.3.1} and the imbedding \eqref{2.3.2}, we have
$$
\Big\|e^{t\Delta}u_0\Big\|_{L^p\big([0, \infty); \dot{H}^{\frac{2+d-p}{p}}_p\big)} 
\simeq \big\|u_0\big\|_{\dot{B}_p^{\frac{d}{p}-1, p}}
\lesssim \big\|u_0\big\|_{\dot{H}_q^{\frac{d}{q} - 1}}. \ \ \ \ \ \ \ \ \ \ \ \ \ \ \ \ \ \  \qed
$$
\begin{Th}\label{Th.5}
Let  $d \geq 3\ and\ 2 < q \leq  d$. Then for any $p$ be such that 
$$
q < p < {\rm min}\Big\{\frac{(d-2) q}{d-q}, \ d+ 2\Big\},
$$
there exists a  constant $\delta_{q, p, d} > 0$ such that for 
all $T > 0$  and for all \linebreak $u_0 \in \dot{H}^{d/q -1}_q(\mathbb{R}^d)$ 
with ${\rm div}(u_0) = 0$  satisfying 
\begin{equation}\label{2.3.3}
\big\|e^{t\Delta}u_0\big\|_{L^p\big([0, T]; \dot{H}^{\frac{2+d-p}{p}}_p\big)}
 \leq \delta_{q, p, d},
\end{equation}
NSE has a unique mild solution
$u \in L^p\big([0, T]; \dot{H}^{\frac{2+d-p}{p}}_p\big) \cap C\big([0, T];
 \dot{H}^{d/q - 1}_q\big)$. Denoting $w = u - e^{t\Delta}u_0$, then we have
$$
w \in L^p\big([0, T]; \dot{H}^{\frac{2+d-p}{p}}_p\big) 
\cap L^\infty\Big([0, T]; \dot{B}^{\frac{d+p-2}{p}-1, 
\frac{p}{2}}_{\frac{d p}{d+p-2}}\Big).
$$
Finally, we have 
\begin{gather*}
\Big\|e^{t\Delta}u_0\Big\|_{L^p\big([0, T]; \dot{H}^{\frac{2+d-p}{p}}_p\big)} 
\leq \big\|u_0\big\|_{\dot{B}_p^{\frac{d}{p}-1, p}} 
\lesssim \big\|u_0\big\|_{\dot{H}_q^{\frac{d}{q} - 1}} , \label{2.3.4}
\end{gather*}
in particular,  for arbitrary $u_0 \in \dot{H}^{d/q -1}_q$ 
the inequality \eqref{2.3.3} holds when $T(u_0)$ is small enough; 
and there exists a positive constant $\sigma_{q,p,d}$ 
such that for all $\big\|u_0\big\|_{\dot{B}_p^{\frac{d}{p}-1, p}}  \leq \sigma_{q,p,d}$
we can take $T = \infty$.
\end{Th}
\textbf{Proof}: The proof of Theorem \ref{Th.5} is similar 
to the one of Theorem  \ref{Th.1}, by combining Theorem \ref{Th.2.1} with 
Lemmas \ref{lem.2.2.6}, \ref{lem.2.1.2} (for $r = p, \ s = \frac{2+d-p}{p}$), 
\ref{lem.2.3.2}, \linebreak and \ref{lem.2.3.3}.\qed
\begin{Not}
The case $q = d$ was treated by several authors, 
see for example (\cite{M. Cannone 1995}, \cite{Hongjie Dong 2007},
\cite{T. Kato 1984}). However their 
results are different from ours. 
\end{Not}
\vskip 0.2cm
\textbf{3.3.2. Solutions to the Navier-Stokes equations 
with initial value in the critical spaces} 
$\dot{H}^{\frac{d}{q} - 1}_q(\mathbb{R}^d)$ 
\textbf{for} $d \geq 3\ {\rm \textbf{and}}\ 1 <  q \leq 2$.
\vskip 0.2cm
\begin{Lem}\label{lem.2.4.2}
Let \ $d \geq 3\ and\ 1 <  q \leq 2$. Then the bilinear 
operator $B(u, v)(t)$ is continuous from 
$$
L^{2q}\Big([0, T]; \dot{H}^{\frac{d + 2 - 2q}{q} }_{\frac{dq}{d + 1 - q}}\Big) 
\times L^{2q}\Big([0, T]; 
\dot{H}^{\frac{d + 2 - 2q}{q} }_{\frac{dq}{d + 1 - q}}\Big) 
$$
into
$$
L^\infty\big([0, T]; \dot{B}^{\frac{d}{q} - 1, q}_q\big),
$$
and we have the inequality
\begin{gather*}
\big\|B(u, v)\big\|_{L^\infty\big([0, T]; \dot{H}^{\frac{d}{q} - 1}_q\big)} 
\lesssim \big\|B(u, v)\big\|_{L^\infty\big([0, T]; 
\dot{B}^{\frac{d}{q} - 1, q}_q\big)} \notag\\
\leq C\Big\|u\Big\|_{L^{2q}\big([0, T]; 
\dot{H}^{\frac{d + 2 - 2q}{q} }_{\frac{dq}{d + 1 - q}}\big)}
\Big\|v\Big\|_{L^{2q}\big([0, T]; 
\dot{H}^{\frac{d + 2 - 2q}{q} }_{\frac{dq}{d + 1 - q}}\big)},
\end{gather*}
where C is a positive constant independent of T.
\end{Lem}
\textbf{Proof}:  Applying Lemma \ref{lem.2.1.3} 
with $r = 2q, \ p = \frac{dq}{d + 1 - q},\ {\rm and}
\ s = \frac{d + 2 - 2q}{q} $, we get
$$
\frac{1}{\tilde{p}} = \frac{2}{p} - \frac{s}{d} =\frac{1}{q},
$$
and from (a) of Lemma \ref{lem.2.2.3}, we have
\begin{gather*}
\big\|B(u, v)\big\|_{L^\infty\big([0, T]; \dot{H}^{\frac{d}{q} - 1}_q\big)} 
\lesssim \big\|B(u, v)\big\|_{L^\infty\big([0, T]; 
\dot{B}^{\frac{d}{q} - 1, q}_q\big)} \notag\\
\lesssim \Big\|u\Big\|_{L^{2q}\big([0, T]; 
\dot{H}^{\frac{d + 2 - 2q}{q} }_{\frac{dq}{d + 1 - q}}\big)}
\Big\|v\Big\|_{L^{2q}\big([0, T]; 
\dot{H}^{\frac{d + 2 - 2q}{q} }_{\frac{dq}{d + 1 - q}}\big)}. 
\ \ \ \ \ \ \ \ \ \ \ \ \ \ \ \ \ \ \ \qed
\end{gather*}
\begin{Lem}\label{lem.2.4.3}
Assume that $u_0 \in \dot{H}^{\frac{d}{q} - 1}_q$ 
with $d \geq 3$ and $1 < q \leq 2$. Then 
$$
\Big\|e^{t\Delta}u_0\Big\|_{L^{2q}\big([0, \infty); 
\dot{H}^{\frac{d + 2 - 2q}{q} }_{\frac{dq}{d + 1 - q}}\big)} 
 \simeq \big\|u_0\big\|_{\dot{B}_{dq/(d + 1 - q)}^{(d + 1)/q - 2, 2q}}
\lesssim \big\|u_0\big\|_{\dot{H}^{d/q - 1}_q}.
$$
\end{Lem}
\textbf{Proof}: By using (a), (c), and (d) of 
Lemma \ref{lem.2.2.3} in order to obtain  
\begin{equation}\label{2.4.2}
 \dot{H}^{\frac{d}{q} - 1}_q \hookrightarrow \dot{B}_q^{\frac{d}{q} - 1, 2}
 \hookrightarrow \dot{B}_q^{\frac{d}{q} - 1, 2q} 
\hookrightarrow \dot{B}_{dq/(d + 1 - q)}^{(d + 1)/q - 2, 2q}.
\end{equation}
Applying Lemma \ref{lem.2.1.4} and from the 
imbedding \eqref{2.4.2} we have the estimates
\begin{gather*}
\Big\|e^{t\Delta}u_0\Big\|_{L^{2q}\big([0, \infty); 
\dot{H}^{\frac{d + 2 - 2q}{q} }_{\frac{dq}{d + 1 - q}}\big)} 
\simeq \big\|\dot{\Lambda}^{\frac{d + 2 - 2q}{q} }
u_0\big\|_{\dot{B}_{dq/(d + 1 - q)}^{-1/q, 2q}} \\
\simeq \big\|u_0\big\|_{\dot{B}_{dq/(d + 1 - q)}^{(d + 1)/q - 2, 2q}} 
\lesssim \big\|u_0\big\|_{\dot{H}_q^{d/q - 1}}. \qed \notag
\end{gather*}
\begin{Th}\label{Th.7}
Let $d \geq 3\ and\ 1 < q \leq 2$. 
There exists a positive constant $\delta_{q, d}$ 
such that for all $T > 0$ and for all 
$u_0 \in \dot{H}^{d/q -1}_q(\mathbb{R}^d)$ 
with ${\rm div}(u_0) = 0$ satisfying 
\begin{equation}\label{2.4.3}
\Big\|e^{t\Delta}u_0\Big\|_{L^{2q}\big([0, T]; 
\dot{H}^{\frac{d + 2 - 2q}{q} }_{\frac{dq}{d + 1 - q}}\big)} \leq \delta_{q, d},
\end{equation}
NSE has a unique mild solution
$u \in L^{2q}\big([0, T]; \dot{H}^{\frac{d + 2 - 2q}{q} }_{\frac{dq}{d + 1 - q}}\big) 
\cap C\big([0, T]; 
\dot{H}^{d/q - 1}_q\big)$. Denoting $w = u - e^{t\Delta}u_0$, then we have
$$
w \in L^{2q}\big([0, T]; \dot{H}^{\frac{d + 2 - 2q}{q} }_{\frac{dq}{d + 1 - q}}\big) 
\cap L^\infty\big([0, T]; \dot{B}^{\frac{d}{q} - 1, q}_q\big).
$$
Finally, we have 
\begin{gather*}
\Big\|e^{t\Delta}u_0\Big\|_{L^{2q}\big([0, T]; 
\dot{H}^{\frac{d + 2 - 2q}{q} }_{\frac{dq}{d + 1 - q}}\big)} 
\leq \big\|u_0\big\|_{\dot{B}_{dq/(d + 1 - q)}^{(d + 1)/q - 2, 2q}} 
\lesssim \big\|u_0\big\|_{\dot{H}^{d/q -1}_q},\notag
\end{gather*}
in particular,  for 
arbitrary $u_0 \in \dot{H}^{d/q -1}_q(\mathbb{R}^d)$ 
the inequality \eqref{2.4.3} holds when $T(u_0)$ is small enough;  
and there exists a positive constant $\sigma_{q, d}$ 
such that for all $\|u_0\|_{\dot{B}_{dq/(d + 1 - q)}^{(d + 1)/q - 2, 2q}}  \leq \sigma_{q, d}$
we can take $T = \infty$.
\end{Th}
\textbf{Proof}: The proof of Theorem \ref{Th.7} is similar to 
the one of Theorem  \ref{Th.1}, by combining Theorem 
\ref{Th.2.1} with Lemmas \ref{lem.2.2.6}, \ref{lem.2.1.2} 
(for $r = 2q, p = \frac{dq}{d+1-q}, s = \frac{d + 2 - 2q}{q}$), 
\ref{lem.2.4.2}, and \ref{lem.2.4.3}. \qed 
\begin{Not}
The case $q = 2$ was treated by several authors, see for example 
(\cite{M. Cannone 1995},\cite{H. Fujita 1964}, \cite{P. G. Lemarie-Rieusset 2002}). 
However their results are different from ours.
\end{Not}
\vskip 0.5cm 
{\bf Acknowledgments}. This research is funded by 
Vietnam National \linebreak Foundation for Science and Technology 
Development (NAFOSTED) under grant number  101.02-2014.50.

\end{document}